\input amstex
\magnification=1200
\documentstyle{amsppt}
\NoRunningHeads
\NoBlackBoxes
\topmatter
\title Games, predictions, interactivity
\endtitle
\author Denis V. Juriev\endauthor
\affil ul.Miklukho-Maklaya 20-180, Moscow 117437 Russia\linebreak
(e-mail: denis\@juriev.msk.ru)
\endaffil
\date math.OC/9906107
\enddate
\keywords Differential games, Interactive games
\endkeywords
\subjclass 90D25 (Primary) 93B52 (Secondary)
\endsubjclass
\abstract\nofrills This short note is devoted to the unraveling of 
the hidden interactivity of ordinary differential games which is 
an artefact of predictions of the behaviour of other players by 
the fixed player and describes deviations of their real behaviour 
from such predictions. A method to improve the predictions is proposed.
Applications to the strategical analysis of interactive games are also 
specified.
\endabstract
\endtopmatter
\document
The mathematical formalism of interactive games, which extends one of ordinary
games (see e.g.[1]) and is based on the concept of an interactive control, 
was recently proposed by the author [2] to take into account the complex 
composition of controls of a real human person, which are often complicated 
couplings of his/her cognitive and known controls with the unknown subconscious 
behavioral reactions. This formalism is applicable also to the description of 
external unknown influences and, thus, is useful for problems in computer 
science (e.g. the semi-artificially controlled distribution of resources), 
mathematical economics (e.g. the financial games with unknown dynamical 
factors) and sociology (e.g. the collective decision making).

However, the interactivity may be unraveled in all differential games. In 
some sense any prediction of behaviour of other players by the fixed player 
allows to consider their controls as interactive. The goal of this article is 
to explain how it is done.

\head I. Interactive systems and games\endhead

\definition{Definition 1 [2]} An {\it interactive system\/} (with $n$
{\it interactive controls\/}) is a control system with $n$ independent 
controls coupled with unknown or incompletely known feedbacks (the feedbacks
as well as their couplings with controls are of a so complicated nature that 
their can not be described completely). An {\it interactive game\/} is a game 
with interactive controls of each player.
\enddefinition

Below we shall consider only deterministic and differential interactive
systems. In this case the general interactive system may be written in the 
form:
$$\dot\varphi=\Phi(\varphi,u_1,u_2,\ldots,u_n),$$
where $\varphi$ characterizes the state of the system and $u_i$ are
the interactive controls:
$$u_i(t)=u_i(u_i^\circ(t),\left.[\varphi(\tau)]\right|_{\tau\leqslant t}),$$
i.e. the independent controls $u_i^\circ(t)$ coupled with the feedbacks on
$\left.[\varphi(\tau)]\right|_{\tau\leqslant t}$. One may suppose that the
feedbacks are integrodifferential on $t$.

However, it is reasonable to consider the {\it differential interactive
games}, whose feedbacks are purely differential. It means that
$$u_i(t)=u_i(u_i^\circ(t),\varphi(t),\ldots,\varphi^{(k)}(t)).$$
A reduction of general interactive games to the differential ones via the 
introducing of the so-called {\it intention fields\/} was described in [2]. 
Below we shall consider the differential interactive games only if the opposite 
is not specified explicitely.

The interactive games introduced above may be generalized in the following 
ways. 

The first way, which leads to the {\it indeterminate interactive games},
is based on the idea that the pure controls $u_i^\circ(t)$ and the 
interactive controls $u_i(t)$ should not be obligatory related in the
considered way. More generally one should only postulate that there are
some time-independent quantities $F_\alpha(u_i(t),u_i^\circ(t),\varphi(t),
\ldots,\varphi^{(k)}(t))$ for the independent magnitudes $u_i(t)$ and 
$u_i^\circ(t)$. Such claim is evidently weaker than one of Def.1. For 
instance, one may consider the inverse dependence of the pure and 
interactive controls: $u_i^\circ(t)=u_i^\circ(u_i(t),\varphi(t),\ldots,
\varphi^{(k)}(t))$.

The inverse dependence of the pure and interactive controls has a nice
psychological interpretation. Instead of thinking of our action consisting
of the conscious and unconscious parts and interpreting the least as 
unknown feedbacks which ``dress'' the first, one is able to consider
our action as a single whole whereas the act of consciousness is in
the extraction of a part which it declares as its ``property''. So
interpreted the function $u_i^\circ(u_i,\varphi,\ldots,\varphi^{(k)})$ 
realizes the ``filtering'' procedure.

The second way, which leads to the {\it coalition interactive games}, is
based on the idea to consider the games with coalitions of actions and to
claim that the interactive controls belong to such coalitions. In this case
the evolution equations have the form
$$\dot\varphi=\Phi(\varphi,v_1,\ldots,v_m),$$
where $v_i$ is the interactive control of the $i$-th coalition. If the 
$i$-th coalition is defined by the subset $I_i$ of all players then
$$v_i=v_i(\varphi(t),\ldots,\varphi^{(k)}(t),u^\circ_j| j\in I_i).$$
Certainly, the intersections of different sets $I_i$ may be non-empty so
that any player may belong to several coalitions of actions. Def.1 gives the
particular case when $I_i=\{i\}$.

The coalition interactive games may be an effective tool for an analysis of
the collective decision making in the real coalition games that spread the
applicability of the elaborating interactive game theory to the diverse 
problems of sociology. 

\remark{Remark 1} One is able to consider interactive games of discrete time 
in the similar manner.
\endremark

\remark{Remark 2} In the most cases one may exclude the first derivative of
$\varphi$ from feedbacks. First, one should transform the feedbacks to
the inverse form, i.e. to express pure controls via interactive ones.
Second, one should substitute the first derivative of $\varphi$ by its
value specified by the evolution equations. Third, one should return to
the direct dependence of controls. To exclude the higher derivatives
it is necessary to perform such procedure several times. However, the
time derivatives of controls will appear. The highest derivative should
be considered as a new control variable as it is often done in the control
theory whereas other derivatives and the control itself will be defined as 
states to make the procedure consistent. Nevertheless, sometimes there some 
difficulties to manipulate with feedbacks with derivatives of $\varphi$. 
One may either practically postulate that the feedbacks depend on $\varphi$ 
only or to consider the discrete time approximation and to use the left 
difference of $\varphi$ in feedbacks and the right difference of $\varphi$ 
in the evolution equations. 
\endremark

\head II. The $\varepsilon$-representations\endhead

Interactive games are games with incomplete information by their nature.
However, this incompleteness is in the unknown feedbacks, not in the 
unknown states. The least situation is quite familiar to specialists in
game theory and there is a lot of methods to have deal with it. For
instance, the unknown states are interpreted as independent controls of
the virtual players and some muppositions on their strategies are done.
To transform interactive games into the games with an incomplete information
on the states one can use the following trick, which is called the 
$\varepsilon$-representation of the interactive game.

\definition{Definition 2} The $\varepsilon$-representation of the 
differential interactive game is a representation of the interactive controls 
$u_i(t)$ in the form
$$u_i(t)=u_i(u^\circ_i(t),\varphi(t),\ldots\varphi^{(k)}(t);\varepsilon_i(t))$$
with the {\sl known\/} function $u_i$ of its arguments $u_i^\circ$,
$\varphi,\ldots,\varphi^{(k)}$ and $\varepsilon_i$, whereas 
$$\varepsilon_i(t)=\varepsilon_i(u^\circ_i(t),\varphi(t),\ldots,\varphi^{(k)}
(t)$$ 
is the {\sl unknown\/} function of $u_i^\circ$ and $\varphi,\ldots,
\varphi^{(k)}$.
\enddefinition

\remark{Remark 3} The derivatives of $\varphi$ may be excluded from the
feedbacks in the way described above. 
\endremark

\remark{Remark 4} One is able to consider the $\varepsilon$-representations
of the indeterminate and coalition interactive games.
\endremark

$\varepsilon_i$ are interpreted as parameters of feedbacks and, thus, 
characterize the internal {\sl existential\/} states of players. It motivates 
the notation $\varepsilon$. Certainly, $\varepsilon$-parameters are not
really states being the unknown functions of the states and pure controls,
however, one may sometimes to apply the standard procedures of the theory
of games with incomplete information on the states. For instance, it is
possible to regard $\varepsilon_i$ as controls of the virtual players.
The na{\"\i}vely introduced virtual players only double the ensemble of 
the real ones in the interactive games but in the coalition interactive 
games the collective virtual players are observed. More sophisticated
procedures generate ensembles of virtual players of diverse structure.

Precisely, if the derivatives of $\varphi$ are excluded from the feedbacks (at 
least, from the interactive controls $u_i$ as functions of the pure controls,
states and the $\varepsilon$-parameters) the evolution equation will have 
the form
$$\dot\varphi(t)=\Phi(\varphi,u_1(u^\circ_1(t),\varphi(t);\varepsilon_1(t)),
\ldots,u_n(u_n^\circ(t),\varphi(t);\varepsilon_n(t))),$$
so it is consistent to regard the equations as ones of the controlled
system with the ordinary controls $u_1,\ldots,u_n,\varepsilon_1,\ldots,
\varepsilon_n$. One may consider a new game postulating that these
controls are independent. Such game will be called the ordinary differential
game associated with the $\varepsilon$-representation of the interactive
game.

\head III. Unraveling the interactivity of ordinary differential 
games\endhead 

Let us consider an arbitrary ordinary differential game with
the evolution equations
$$\dot\varphi=\Phi(\varphi,u_1,u_2,\ldots,u_n),$$
where $\varphi$ characterizes the state of the system and $u_i$ are
the ordinary controls. Let us fix a player. For simplicity of notations we
shall suppose that it is the first one. As a rule the players have their
algorithms of predictions of the behaviour of other players. For a fixed
moment $t_0$ of time let us consider the prediction of the first player for
the game. It consists of the predicted controls $u^\circ_{[t_0];i}(t)$ 
($t>t_0$; $i\ge2$) of all players and the predicted evolution of the system 
$\varphi^\circ_{[t_0]}(t)$. Let us fix $\Delta t$ and consider the real and 
predicted controls for the moment $t_0+\Delta t$. Of course, they may be
different because other players use another algorithms for the game prediction.
One may interpret the real controls $u_i(t)$ ($t=t_0+\Delta t$; $i\ge2$) of 
other players as interactive ones whereas the predicted controls 
$u^\circ_{[t_0];i}(t)$ as pure ones, i.e. to postulate their relation in the 
form:
$$u_i(t)=u_i(u^\circ_{[t-\Delta t];i}(t);\varphi^\circ_{[t_0]}(\tau)|\tau\le t).$$
In particular, the feedbacks may be either reduced to differential form via
the introducing of the intention fields or simply postulated to be differential.
Thus, we constructed an interactive game from the initial ordinary game.
One may use $\varphi(\tau)$ as well as $\varphi^\circ_{[t_0]}(\tau)$ 
in the feedbacks. 

Note that the controls of the first player may be also considered as
interactive if the corrections to the predictions are taken into account
when controls are chosen.
              
The obtained construction may be used in practice to make more adequate
predictions. Namely, {\sl a posteriori\/} analysis of the differential
interactive games allows to make the short-term predictions in such games
[3]. One should use such predictions instead of the initial ones. Note that
at the moment $t_0$ the first player knows the pure controls of other
players at the interval $[t_0,t_0+\Delta t]$ whereas their real freedom
is interpreted as an interactivity of their controls $u_i(t)$. So it is
reasonable to choose $\Delta t$ not greater than the admissible time depth 
of the short-term predictions. Estimations for this depth were proposed
in [3]. 

Na{\"\i}vely, the proposed idea to improve the predictions is to consider
deviations of the real behaviour of players from the predicted ones as
a result of the interactivity, then to make the short-term predictions
taking the interactivity into account and, thus, to receive the 
corrections to the initial predictions. Such corrections may be regarded
as ``psychological'' though really they are a result of different methods
of predictions used by players.

Such procedure can be performed also for the discrete time games and, 
thus, is applicable to a wide class of the model entertainment games such 
as domino or chess. On the other hand, the procedure is useful for almost 
all practically important games.

Constructions above motivates the following definition.

\definition{Definition 3} The {\it scenario\/} of an interactive game is
the set of pure controls $u_i^\circ(t)$ as functions of time. The {\it
performance\/} is an interactive game with {\sl a priori\/} fixed scenario.
\enddefinition

Thus, one may said that at any moment $t_0$ an ordinary differential game 
coincides with certain performance in the nearest future $t\in[t_0,t_0+
\Delta t]$.

\remark{Remark 5}
Performances have some specific features. The real evolution of the game
does not {\sl strategically\/} differ from the scenario. Their divergence
is only {\sl tactical\/} that strengthens the role of short-term predictions
in the analysis of game. It explains why ordinary games can be represented
as performances only in the nearest future.
\endremark

\remark{Remark 6}
The interpretation of the ordinary differential game as an interactive game
also allows to perform the strategical analysis of interactive games. Indeed,
let us consider an arbitrary differential interactive game $A$. Specifying its
$\varepsilon$-representation one is able to construct the associated 
ordinary differential game $B$ with the doubled number of players. Making
some predictions in the game $B$ one transform it back into an interactive 
game $C$. Combination of the strategical long-term predictions in the game 
$B$ with the short-term predictions in $C$ is often sufficient to obtain the 
adequate strategical prognosis for $A$.
\endremark

\remark{Remark 7} Some procedures of predictions are applicable only to
the proper classes $U$ of scenarios and specific initial positions. For 
instance, one may consider such class of scenarios and initial positions 
for which certain quantities $Z(t)=Z(u^\circ_i(t),\mathbreak\varphi(t),
\varepsilon_i(t))$ are time-independent or their time derivatives are 
expressed via $Z$ themselves. Besides local-in-time expressions the functional 
expressions of the form $Z(t)=Z(u^\circ_i(t\!-\!\tau),\varphi(t\!-\!\tau),
\varepsilon_i(t\!-\!\tau)|\tau\!\in\![0,\Delta t])$ are also admitted. Also 
one may use the generalized functions (distributions) of the same arguments. 
The knowledge of such quantities allows to make some specific predictions and 
their detection in the concrete game may be rather natural. For instance, in 
the kaleidoscope-roulettes [4] the appearing of similar quantities (but of 
discrete time) indicates the presence of resonances for the class of scenarios 
that allows to make predictions, on the other hand, they also describe the 
interpretational figures detected visually by the player, which should be 
regarded therefore as {\sl ``omens''}. The unraveling of such {\sl omens\/} and their 
interpretation provide a way to success in the kaleidoscope-roulettes. Because 
the controls of various players are independent it is important to know whether 
the situation is stable, it means that if one consider a scenario 
$\{u_i^\circ(t)\}$ in the neigbourhood of the class $U$ (in certain topology 
on the space of all scenario) then slightly improving $Z$ it is possible to 
provide its invariance for this scenario. Bifurcations and their controlling 
are also interesting. For instance, the controlling of bifurcations of $Z$ 
may be interpreted as the {\sl ``omen-manipulation''\/} in the 
kaleidoscope-roulettes that is perhaps the main intrigue in these games.
\endremark

\head IV. Conclusions\endhead

Thus, the hidden interactivity of the ordinary differential games is 
unraveled. It is shown that such games may be regarded as interactive
games of special kind. In these interactive games the pure controls
of players coincide with the predicted ones and, hence, are known for 
the nearest future $t\in[t_0,t_0+\Delta t]$, whereas the real freedom 
of players is interpreted as an interactivity of their really observed 
controls. Such observation leads to a method which allows to improve 
the initial predictions for the game. Applications to the strategical
analysis of interactive games are also briefly specified.

\Refs
\roster
\item"[1]" Isaaks R., Differential games. Wiley, New York, 1965;\newline
Owen G., Game theory, Saunders, Philadelphia, 1968;\newline
Vorob'ev N.N., Current state of the game theory, Uspekhi Matem. Nauk 15(2) 
(1970) 80-140 [in Russian].
\item"[2]" Juriev D., Interactive games and representation theory. I,II.
E-prints: math.FA/9803020, math.RT/9808098.
\item"[3]" Juriev D., The laced interactive games and their {\sl a
posteriori\/} analysis. E-print:\linebreak math.OC/9901043; Differential
interactive games: The short-term predictions. E-print: math.OC/9901074.
\item"[4]" Juriev D., Perceptive games, the image understanding and 
interpretational geometry. E-print: math.OC/9905165; Kaleidoscope-roulette:
the resonance phenomena in perception games. E-print: math.OC/9905180.
\endroster
\endRefs
\enddocument